\RequirePackage{ifpdf}
\ifpdf 
\documentclass[pdftex]{sigma}
\else
\documentclass{sigma}
\fi

\usepackage{euscript}

\newlength{\setBracketHeight}

\newcommand{\LieDer}{\ensuremath{\EuScript L}}
\newcommand{\hook}{\ensuremath{\mathbin{ \hbox{\vrule height1.4pt
        width4pt depth-1pt \vrule height4pt width0.4pt depth-1pt}}}}

\newcommand{\C}[1]{\ensuremath{\mathbb{C}^{#1}}}

\newcommand{\GL}[1]{\operatorname{GL}\left({#1}\right)}

\newcommand{\Lm}[2]{\ensuremath{\Lambda^{#1} \left ( {#2} \right )}}
\newcommand{\nForms}[2]{\ensuremath{\Omega^{#1} \left ( {#2} \right
    )}}

\DeclareMathOperator{\Ad}{Ad}

\DeclareMathOperator{\Aut}{Aut}

\newcommand{\Proj}[1]{\mathbb{P}^{#1}}

\newcommand{\OO}[1]{
  \ensuremath{
    \mathcal{O}
    \ifthenelse{\equal{#1}{0}}
      {}
      {\left({#1}\right)}
  }
}
\newcommand{\OOp}[2]{
  \ensuremath{
    \mathcal{O}
    \ifthenelse{\equal{#1}{0}}
      {}
      {\left({#1}\right)}
    \ifthenelse{\equal{#2}{1}}
      {}
      {^{\oplus{#2}}}
  }
}

\newcommand{\extraSection}[2]
{
\ifthenelse{\boolean{abridged}}
  {
  }
  {
    \section{#1}
    \begin{center}
    \emph{This section will not be referred to
    subsequently, and may be skipped.}
    \end{center}
    \par{#2}
  }
}
\newcommand{\extraSubsection}[2]
{
\ifthenelse{\boolean{abridged}}
  {
  }
  {
    \subsection{#1}
    \begin{center}
    \emph{This subsection will not be referred to
    subsequently, and may be skipped.}
    \end{center}
    \par{#2}
  }
}
\newcommand{\extraStuff}[1]
{
\ifthenelse{\boolean{abridged}}
  {
  }
  {
    {#1}
  }
}

\newcommand{\LieG}{\ensuremath{\mathfrak{g}}}
\newcommand{\LieP}{\ensuremath{\mathfrak{p}}}
\newcommand{\LieM}{\ensuremath{\mathfrak{m}}}
\newcommand{\LieA}{\ensuremath{\mathfrak{a}}}
\newcommand{\LieN}{\ensuremath{\mathfrak{n}}}
\newcommand{\LieH}{\ensuremath{\mathfrak{h}}}
\newcommand{\LieL}{\ensuremath{\mathfrak{l}}}

\newcommand{\Gtot}{\ensuremath{G}}
\newcommand{\gtot}{\mathfrak{g}}
\newcommand{\Gstruc}{P}
\newcommand{\gstruc}{\mathfrak{p}}

\newcommand{\RedPart}{LA}
\newcommand{\redpart}{\mathfrak{l} \oplus \mathfrak{a}}

\newcommand{\cartan}{\ensuremath{\mathfrak{h}}}

\newcommand{\cpt}[1]{%
        \ifthenelse{\equal{#1}{\RedPart} \or \equal{#1}{\redpart}}
        {\ensuremath{(#1)^c}}
        {\ensuremath{#1^c}}
}

\begin{document}

\allowdisplaybreaks

\renewcommand{\PaperNumber}{090}

\FirstPageHeading

\ShortArticleName{Holomorphic Parabolic Geometries and Calabi--Yau Manifolds}

\ArticleName{Holomorphic Parabolic Geometries\\ and Calabi--Yau Manifolds}

\Author{Benjamin MCKAY}

\AuthorNameForHeading{B.~McKay}

\Address{School of Mathematical Sciences, University College Cork, Cork, Ireland}
\Email{\href{mailto:b.mckay@ucc.ie}{b.mckay@ucc.ie}}
\URLaddress{\url{http://euclid.ucc.ie/pages/staff/Mckay/}}

\ArticleDates{Received May 25, 2011, in f\/inal form September 15, 2011;  Published online September 20, 2011}

\Abstract{We prove that the only complex parabolic geometries on
Calabi--Yau manifolds are the homogeneous geometries on complex
tori. We also classify the complex parabolic geometries on
homogeneous compact K\"ahler manifolds.}


\Classification{53C55; 53A55; 53C10}

\section{Introduction}

We will prove that Calabi--Yau manifolds (other
than those covered by complex tori) cannot bear
holomorphic parabolic geometries.
Gunning \cite{Gunning:1978} proved that any compact K{\"a}hler surface
with $c_1=0$ admitting
a holomorphic normal projective or conformal connection
is covered by a complex torus. Kobayashi \cite{Kobayashi/Horst:1983}
proved that any compact K{\"a}hler manifold
with $c_1=0$ which admits a holomorphic normal projective or conformal connection
is covered by a complex torus. Our arguments are simpler than
those of Gunning \cite{Gunning:1978} or Kobayashi \cite{Kobayashi/Horst:1983},
and give stronger conclusions (not requiring
normalcy, and applying directly to all
holomorphic parabolic geometries, not just projective and
conformal connections).

\section{Review of the literature}

Let us contrast our results in this paper with those
of \cite{Biswas/McKay:2010a,Biswas/McKay:2010b,Dumitrescu:2010}.
In \cite{Biswas/McKay:2010a},
we proved that if a~smooth complex projective variety
with $c_1=0$ bears a holomorphic Cartan geometry, then the
smooth projective variety has holomorphic unramif\/ied covering
map by an Abelian variety. In~\cite{Biswas/McKay:2010b},
we generalized this result to prove that
if a compact K{\"a}hler manifold with $c_1=0$
bears a~holomorphic Cartan geometry, then the
compact K{\"a}hler manifold has a holomorphic unramif\/ied covering
map by a complex torus. A special case of this result
(for parabolic geometries) will be proven
here in the f\/irst part of Theorem~\ref{thm:CY}.
Sorin Dumitrescu proved this same result,
and a collection of remarkable generalizations (for
example obstructions to holomorphic
Cartan geometries on products~$CY \times X$).
However, those proofs used algebraic geometry
where the proof below relies more on the local theory
of parabolic geometries. Dumitrescu also
proved similar results to those below
on Calabi--Yau manifolds bearing
structures of af\/f\/ine algebraic type~\cite{Dumitrescu:2001}.
The paper you are currently reading improves in one way
on all of those results:  Theorem~\ref{theorem:translation}
shows that holomorphic parabolic geometries
on tori are translation invariant. Moreover the proof
of Theorem~\ref{theorem:translation} uses exactly
the same local calculations as the proof of Theorem~\ref{thm:CY}.
It is well known
that there are holomorphic Cartan geometries on complex
tori which are not translation invariant, so Theorem~\ref{theorem:translation}
is surprising. The
classif\/ication of holomorphic parabolic geometries
on compact K{\"a}hler manifolds with $c_1=0$
is completed below (up to f\/inite unramif\/ied covering),
in Theorems~\ref{thm:CY} and~\ref{theorem:translation}.
This classif\/ication is then use to classify
holomorphic parabolic geometries on homogeneous
compact K{\"a}hler manifolds.
The analogous classif\/ications for holomorphic Cartan
geometries are not known or conjectured.

Let us review the current state of the search
for holomorphic Cartan geometries on compact
complex manifolds. Holomorphic Cartan geometries
are completely classif\/ied on (1)~any compact Riemann
surface \cite{McKay:2011}, (2)~any compact
complex surface containing a rational curve~\cite{McKay:2011}
and (3)~any
compact rationally connected complex manifold~\cite{Biswas/McKay:2010c}. (This last class of manifolds
includes, for example, all Fano manifolds, and all
rational homogeneous varieties $G/P$.)
On compact K{\"a}hler manifolds with $c_1=0$, we classify
below all of the holomorphic parabolic geometries.
Suppose that $P \subset G$ is a maximal parabolic subgroup
of a complex semisimple Lie group.
For any compact K{\"a}hler manifold~$M$ containing a rational
curve, either $M=G/P$ or else~$M$ admits no holomorphic $G/P$-geometry~\cite{Biswas/McKay:2010c}.

On every compact complex surface, all holomorphic torsion-free
af\/f\/ine connections are classif\/ied~\cite{Klingler:1998,Dumitrescu:2010b}.
On every compact complex surface,
all holomorphic normal projective connections
which do not arise from holomorphic
af\/f\/ine connections are classif\/ied~\cite{Klingler:1998,Dumitrescu:2010b}.
(It remains to see which pairs of holomorphic af\/f\/ine connections
determine the same holomorphic projective connection.)
It is known which compact complex surfaces admit
holomorphic normal parabolic geometries
\cite{Inoue/Kobayashi/Ochiai:1980,Kobayashi/Horst:1983,%
KobayashiOchiai:1980,Kobayashi/Ochiai:1981ii,Kobayashi/Ochiai:1982}.

Suppose that $M$ is a locally symmetric
complex manifold of f\/inite volume,
$M=\Gamma \backslash X$, where $X$
is a noncompact Hermitian symmetric space, say with
compact dual $G/P$,
and $\Gamma$ is a discrete group of isomorphisms
of the standard f\/lat $G/P$-geometry on~$X$.
Suppose further
that $G/P$ is an irreducible symmetric space.
Then $M$ admits a unique normal holomorphic
parabolic geometry modelled on $G/P$ (the obvious one)
\cite{Kobayashi/Ochiai:1981,Klingler:2001}.

On any compact K{\"a}hler manifold,
there are constraints on characteristic classes
arising from the presence of holomorphic Cartan
geometries \cite{Kobayashi/Ochiai:1981,McKay:2011c}.
The smooth complex projective 3-folds that bear
holomorphic normal projective or
conformal connections are classif\/ied
\cite{Jahnke/Radloff:2002,Jahnke/Radloff:2004}.

Every known example of a compact complex manifold
admitting a holomorphic Cartan geo\-met\-ry
also admits a f\/lat holomorphic Cartan geometry
with the same model with one exception: there are translation
invariant Cartan geometries on the complex torus (see below),
which are not f\/lat, and for which no complex torus admits
a f\/lat geometry with the same model.
Many compact complex manifolds only admit locally
homogeneous holomorphic Cartan geometries, but some
also admit locally inhomogeneous ones~\cite{Dumitrescu:2010b}.
There are sporadic results classifying f\/lat holomorphic
Cartan geometries of various types on various complex
manifolds \cite{McKay/Pokrovskiy:2010}, but
I am not aware of any other results concerning the
classif\/ication of holomorphic Cartan geometries.

\section[Calabi-Yau manifolds]{Calabi--Yau manifolds}

\begin{definition} For this article, a \emph{Calabi--Yau manifold}
is a compact K{\"a}hler manifold $M$ with $c_1(TM)=0$.
\end{definition}

It is well known that a Calabi--Yau manifold satisf\/ies $c_2(TM)=0$ just
if it has a torus as unramif\/ied covering space.
Let us recall how this follows from Yau's proof of the Calabi conjecture.
For any K{\"a}hler manifold, say of dimension $n$,
with $\Omega$ its K{\"a}hler form, it is easy to calculate that
\[
c_2 \wedge \Omega^{n-2} = \big(\left\|R\right\|^2 + \text{scalar}^2
- 2 \left\|\text{Ricci}\right\|^2 \big) \Omega^n,
\]
(see Berger and Lascoux \cite{Lascoux/Berger:1970})
where $R$ is the curvature tensor. If $c_1=0$, then
there is a metric for which $\text{Ricci}=0$, by
Yau's solution of the Calabi conjecture~\cite{Yau:1977}.
Hence $c_2=0$ implies $R=0$, f\/lat. But then
$M$ is covered by a f\/lat
torus (see Igusa~\cite{Igusa:1954}).

\begin{lemma}[Inoue, Kobayashi and Ochiai \cite{Inoue/Kobayashi/Ochiai:1980}]
Any compact complex manifold which bears a holomorphic
Cartan geometry with reductive algebraic structure
group has vanishing Atiyah class.
In particular, if K{\"a}hler then
it is the quotient of a complex torus
under a finite unramified covering map.
\end{lemma}
\begin{proof}
The Cartan connection splits invariantly into a sum
of a connection (in the sense of Ehresmann) and
a soldering form; see Sharpe \cite[Lemma~2.1, p.~362]{Sharpe:1997}.
The existence of a~connection is precisely the
vanishing of the Atiyah class; see Atiyah~\cite{Atiyah:1957}.
If K{\"a}hler, then all Chern classes of the
tangent bundle vanish just when the Atiyah class
does. By the previous discussion, the manifold has a torus as f\/inite unramif\/ied
covering space.
\end{proof}

\begin{example}
A holomorphic Riemannian metric is a simple
example of a reductive Cartan geometry,
and our results tell us that holomorphic
Riemannian metrics can not live on any
compact K{\"a}hler manifold except those
covered by tori. This is well known (see
Inoue, Kobayashi and Ochiai~\cite{Inoue/Kobayashi/Ochiai:1980}).
\end{example}

We will prove the following theorem.

\begin{theorem}\label{thm:CY}
If a Calabi--Yau manifold bears a holomorphic parabolic
geometry, then it is covered by a torus. More generally,
any compact complex manifold with
a holomorphic parabolic geometry and trivial
canonical bundle must have
a holomorphic affine connection.
\end{theorem}

\section{Rational homogeneous varieties}

Suppose that $G/P$ is a rational homogeneous variety,
so $G$ is a complex semisimple Lie group
and $P$ is a complex parabolic subgroup,
with Lie algebras $\LieG$ and $\LieP$. We can
express $\LieP$ as a sum
of the Cartan subalgebra of $\LieG$ together with
various root spaces, including all of the positive
root spaces. Some negative root spaces will also lie
in $\LieP$. Once we f\/ix the choice
of $\LieG$ and $\LieP$, roots then divide up into
3 categories as follows. The \emph{compact roots} of $\LieG$ are
the roots $\alpha$ of $\LieG$ so that the root spaces of
both $\alpha$ and $-\alpha$ belong to the Lie algebra
of $\LieP$.
All other roots are \emph{noncompact}, and divide
into the \emph{noncompact positive} and \emph{noncompact negative
roots}, according to whether or not their root spaces
lie in $\LieP$.
The Dynkin diagram of $G/P$ is the Dynkin
diagram of $G$ (labelled by simple roots),
with simple roots dotted if they are compact,
and crossed if they are noncompact.

The sum of the root spaces of the noncompact positive roots
is the maximal nilpotent subalgebra, denoted
$\LieN \subset \LieP$.
The sum of the root spaces of the noncompact negative roots
is also a nilpotent subalgebra, denoted
$\LieN^{-} \subset \LieG$, complementary to $\LieP$.
Let $\LieA \subset \LieP$
be the subalgebra spanned by the coroots of the compact roots.
Let $\LieM \subset \LieP$ be the Lie subalgebra generated
by the root space of the compact roots.
The Lie subalgebra
$\LieM \oplus \LieA$ ($\LieN$ resp.)
is the maximal reductive (nilpotent) subalgebra
of $P$; see Knapp \cite{Knapp:2002}.
Let $M$, $A$, $N$ and $N^-$ be the connected subgroups
of $G$ with Lie algebras
$\LieM$, $\LieA$, $\LieN$ and $\LieN^-$.
The groups $M$, $A$, $N$, $N^-$ and $P$
are all algebraic (see Fulton and Harris \cite[p.~382]{Fulton/Harris:1991}).
The splitting
$\LieG = \LieN^{-} \oplus \LieM \oplus \LieA \oplus \LieN$
is $MA$-invariant.

Pick a Chevalley basis $X_{\alpha}$, $H_{\alpha}$
for $\LieG$. Recall (see Serre~\cite{Serre:2001}) that this is a basis
parameterized by roots $\alpha \in \cartan^*$
(with $\cartan \subset \gtot$ a Cartan
subalgebra) for which
\begin{enumerate}\itemsep=0pt
\item $\left[H,X_{\alpha}\right]=\alpha(H) X_{\alpha}$
for each $H \in \cartan$;
\item $\alpha\left(H_{\beta}\right)=
2 \frac{\left<\alpha, \beta\right>}{\left<\beta,\beta\right>}$
(measuring inner products via the Killing form);
\item $\left[H_{\alpha},H_{\beta}\right]=0$;
\item
\[
\left[X_{\alpha},X_{\beta}\right]=
\begin{cases}
H_{\alpha}, & \text{ if } \alpha+\beta=0, \\
N_{\alpha\beta} X_{\alpha+\beta}, & \text{otherwise}
\end{cases}
\]
with
\begin{enumerate}\itemsep=0pt
\item $N_{\alpha\beta}$ an integer,
\item $N_{-\alpha,-\beta}=-N_{\alpha\beta}$,
\item If $\alpha$, $\beta, $ and $\alpha+\beta$ are roots, then
$N_{\alpha\beta}=\pm (p+1)$,
where $p$ is the largest integer for which $\beta-p  \alpha$ is a root,
\item
$N_{\alpha \beta}=0$ if $\alpha+\beta = 0$ or
if any of $\alpha$, $\beta,$ or  $\alpha+\beta$ is not a root.
\end{enumerate}
\end{enumerate}

Consider the 1-forms $\omega^{\alpha}$ dual to the vectors
$X_{\alpha}$ of a Chevalley basis. We use the Killing form to extend $\alpha$ from $\cartan$ to $\gtot$,
by splitting $\gtot = \cartan + \cartan^{\perp}$,
and taking $\alpha=0$ on $\cartan^{\perp}$. The 1-forms~$\omega^{\alpha}$,~$\alpha$
span $\gtot^*$.

Each exterior form in $\Lm{*}{\gtot}^*$ extends uniquely to
a left invariant dif\/ferential form in $\nForms{*}{\Gtot}$, and we will
identify these. These forms determine a basis of left invariant 1-forms
$\omega^{\alpha}$, $\alpha$, and a~basis of left invariant vector
f\/ields $X_{\alpha}$, $H_{\alpha}$. Clearly
\begin{gather*}
d \omega^{\alpha}  = - \alpha \wedge \omega^{\alpha} - \frac{1}{2}
\sum_{\beta + \gamma=\alpha} N_{\beta \gamma} \omega^{\beta} \wedge \omega^{\gamma}, \qquad
d \alpha  =
-
\sum_{\beta} \frac{\left<\alpha,\beta\right>}
{\left<\beta,\beta\right>} \omega^{\beta} \wedge \omega^{-\beta},
\end{gather*}
with sums over all roots.
To be more precise $\omega^{\alpha}$, $\alpha$ is not quite a basis of 1-forms,
since there will be relations among the $\alpha$ 1-forms
in general. To produce a basis, we would have to restrict
to the $\alpha$ 1-forms which are simple roots, but
include all of the $\omega^{\alpha}$ 1-forms,
even for nonsimple~$\alpha$. The basis~$\omega^{\alpha}$,~$\alpha$
is \emph{not} the dual basis to $X_{\alpha}$, $H_{\alpha}$.

\begin{definition}
If $G/P$ is a rational homogeneous variety,
let $\delta = \delta_{G/P}$ be
\[
\delta = \frac{1}{2} \sum_{\alpha}^{\times} \alpha,
\]
where $\sum^{\times}$ means the sum over all noncompact negative roots.
\end{definition}

\begin{lemma}
The Killing form inner product
\(
\left<\delta, \beta\right>
\)
$($where $\delta$ is half the sum of noncompact negative roots, and $\beta$
any root$)$ vanishes
just precisely for $\beta$ a root of the
maximal semisimple subalgebra $\mathfrak{m} \subset \gstruc$.
\end{lemma}
\begin{proof} Knapp \cite[Corollary 5.100, p.~330]{Knapp:2002}
gives a completely elementary proof. We give a~proof along
the same lines, to keep our exposition self-contained.
Pick $\beta$ any positive root. If~$\gamma$ is any
noncompact negative root, and $\left<\gamma,\beta\right> > 0$,
then
\[
\gamma, \gamma-\beta, \gamma-2   \beta, \dots, \gamma-q   \beta = r_{\beta} \gamma
\]
is a string of roots ending in the ref\/lection $r_{\beta}$ of $\gamma$.
To start with, $\gamma$ already contains a positive
multiple of a noncompact negative simple root. Equivalently, $\gamma$
has some negative multiple of a~noncompact positive simple root~$\alpha_1$.
Subtracting the positive
root~$\beta$ can only make the multiple of~$\alpha_1$
larger negative. Therefore the entire string consists
of noncompact negative roots.

If we have an entire $\beta$-string of noncompact
negative roots, for a positive
root $\beta$, clearly
\begin{gather*}
\left<r_{\beta} \gamma,\beta\right>
 = -\left<r_{\beta} \gamma, r_{\beta} \beta\right>
 =
-\left<\gamma,\beta\right>.
\end{gather*}
Therefore $\left<\gamma,\beta\right>$ cancels with $\left<r_{\beta}
\gamma,\beta\right>$ in the sum $\left<\delta,\beta\right>$. Hence
the entire string cancels out of that sum.

It follows that
\begin{equation}\label{eqn:Sum}
\left<\delta,\beta\right>=
\sum_{\gamma} \left<\gamma,\beta\right>
\end{equation}
where the sum is over noncompact negative roots $\gamma$
for which both $\left<\gamma,\beta\right> \le 0$
and for which the other end of the $\beta$-string
through $\gamma$ is noncompact positive.
Of course, the terms with  $\left<\gamma,\beta\right> = 0$
cancel out too, so the sum~(\ref{eqn:Sum}) is over
noncompact negative roots $\gamma$ for
which both $\left<\gamma,\beta\right> < 0$
and for which the other end of the $\beta$-string
through~$\gamma$ is noncompact positive. In particular,
the sum (\ref{eqn:Sum}) is a sum of negative terms.
But there might not be any terms.

If $\beta$ is a compact root, then clearly ref\/lection in
$\beta$ preserves the roots belonging to the parabolic subalgebra
$\gstruc$, and therefore preserves the noncompact negative roots.
So the noncompact negative roots will all lie in $\beta$-strings, and
$\left<\delta,\beta\right>=0$ for these roots. On the other hand, if~$\beta$ is a noncompact root, then $\beta$
is either positive or negative. We can assume that~$\beta$ is positive, since we only need to show that $\left<\delta,\beta\right> \ne 0$. Take
$\gamma=-\beta$, to see that the sum~(\ref{eqn:Sum}) has at least
one negative term.
\end{proof}

\section{Parabolic geometries}

The standard reference on parabolic
geometries is \v{C}ap and Slovak \cite{Cap/Slovak:2009};
we will use the standard def\/initions, as in their
book, which are far too long to put into this paper.
Suppose that $E \to M$ is a holomorphic
parabolic geometry, with some model $G/P$.
Very similar structure equations hold for any
holomorphic parabolic geometry with the same model. Indeed,
the Cartan connection is a 1-form valued in the Lie
algebra $\mathfrak{g}$ of $G$, so splits into a sum
of 1-forms~$\omega^{\alpha}$ and~$\alpha$ from the
decomposition of $\mathfrak{g}$ into root spaces.
From the def\/inition of a Cartan geometry, the
Cartan connection satisf\/ies the same structure
equations as the Maurer--Cartan form on the model,
but with semibasic curvature correction terms, so
\begin{gather*}
d \omega^{\alpha}  = - \alpha \wedge \omega^{\alpha} - \frac{1}{2}
\sum_{\beta + \gamma=\alpha} N_{\beta \gamma} \omega^{\beta} \wedge \omega^{\gamma}
+
\sum^{\times}_{\beta, \gamma}
\kappa^{\alpha}_{\beta \gamma} \omega^{\beta} \wedge \omega^{\gamma},
\\
d \alpha  =
-
\sum_{\beta} \frac{\left<\alpha,\beta\right>}
{\left<\beta,\beta\right>} \omega^{\beta} \wedge \omega^{-\beta}
+
\sum^{\times}_{\beta, \gamma}
\lambda^{\alpha}_{\beta \gamma} \omega^{\beta} \wedge \omega^{\gamma},
\end{gather*}
where the $\kappa$ and $\lambda$ terms are Cartan geometry
curvature terms, so they vanish except possibly for
$\beta$ and $\gamma$ noncompact negative roots, and
once again $\displaystyle{\sum^{\times}}$ means
the sum over noncompact negative roots.

It is vital in the following that, even if we
work on a manifold where we have imposed some
relations on these 1-forms, we will still use the
Killing form on the original Lie algebra $\gtot$
to compute inner products $\left<\alpha,\beta\right>$.
This is our only notational ambiguity.

\section{Proofs of the theorems}

Replacing our Calabi--Yau manifold by a f\/inite covering space
if needed, we can assume that it
bears a nowhere-vanishing holomorphic volume
form. We then derive our theorem from the following stronger theorem:

\begin{theorem}\label{theorem:Reduction}
If a complex manifold bears a holomorphic parabolic geometry
and a holomorphic volume form, then it admits a canonical
holomorphic reduction of the structure group of the parabolic
geometry to a reductive algebraic group.
\end{theorem}
\begin{proof}
Suppose that $E \to M$ is a holomorphic parabolic
geometry modelled on $\Gtot/\Gstruc$,
and $\sigma$ a~holomorphic volume form on $M$.
Pick a Chevalley basis. Let
\begin{gather*}
\Omega  = \bigwedge^{\times}_{\alpha} \omega^{\alpha}, \qquad
\delta  = \frac{1}{2} \sum^{\times}_{\alpha} \alpha,
\end{gather*}
where the wedge product and sum are over noncompact negative roots.
The sign of $\Omega$ depends on a choice of ordering
of the noncompact negative roots, but any ordering can be chosen, as
long as we are consistent.

We claim that $d \Omega = - 2 \delta \wedge \Omega$.
Order the noncompact negative
roots arbitrarily as $\alpha_1, \alpha_2, \dots$.
Expand out $d \Omega$:
\[
d \Omega
=
 \sum_{j}
(-1)^{j+1}
\bigwedge_{i < j} \omega^{\alpha_i}
\wedge d \omega^{\alpha_j}
\wedge
\bigwedge_{i > j} \omega^{\alpha_i},
\]
by passing the exterior derivative operator along
the various factors, hitting one $\omega^{\alpha}$
at a time, and sticking a suitable $\pm$ sign in front.
Plug in the equation for $d \omega^{\alpha}$,
and the curvature terms all vanish because they
occur in pairs of $\omega^{\beta} \wedge \omega^{\gamma}$,
and at least one of these $\omega^{\beta}$ or $\omega^{\gamma}$
is still present in factors in that term.
If we f\/ind a term like $N_{\beta \gamma} \omega^{\beta}
\wedge \omega^{\gamma}$ in $d \omega^{\alpha}$,
then we must have $\alpha=\beta+\gamma$.
Since $\alpha$ is a noncompact negative root, at least one
of $\beta$ and $\gamma$ must be as well. Neither can be
equal to $\alpha$, since $N_{0\alpha}=0$. Therefore
each such term drops out of $d \Omega$. The reader
now only has to check the signs to see that
$d \Omega = -2 \delta \wedge \Omega$.

Our nonzero section $\sigma$ of the canonical
bundle of $M$ can be pulled back to $E$ as
\(
\sigma = s \Omega,
\)
for a unique nowhere-vanishing function $s : E \to \C{}$.
If $\sigma$ is holomorphic, then
\begin{gather*}
0  =
d\sigma
 =
ds \wedge \Omega + s   d \Omega
=
\left(ds - 2 s \delta \right) \wedge \Omega.
\end{gather*}
Let $P_0 \subset P$ be the subgroup of $P$
acting trivially on
$\Lambda^{\text{top},0}\left(\mathfrak{g}/\mathfrak{p}\right)$.
Let $E' \subset E$ be the set of points at which
$s=1$. Then $E' \subset E$ is a smooth hypersurface since $ds \ne 0$
on tangent spaces of~$E$ along~$E'.$ Clearly
$E'$ is a principal right $P_0$-bundle.

On $E'$, $\delta \wedge \Omega = 0$.
Therefore $\delta$ is semibasic on $E'$:
\[
\delta = \sum^{\times}_{\alpha} t_{\alpha} \omega^{\alpha},
\]
a sum over noncompact negative roots $\alpha$, for some
functions $t_{\alpha} : E' \to \mathbb{C}.$
Taking exterior derivative, we f\/ind
\begin{gather*}
0 =
d \left(\delta - \sum _{\alpha} t_{\alpha} \omega^{\alpha} \right)
 = \sum_{\alpha}
\left(
d \alpha
-
dt_{\alpha} \wedge \omega^{\alpha}
- t_{\alpha} d \omega^{\alpha}
\right)
\\
\phantom{0}{} =
-
\sum_{\beta} \frac{\left<\delta,\beta\right>}
{\left<\beta,\beta\right>} \omega^{\beta} \wedge \omega^{-\beta}
-
\sum_{\alpha}
\left(dt_{\alpha} - t_{\alpha} \alpha\right) \wedge \omega^{\alpha}
\\
\phantom{0=}{} - \frac{1}{2}
\sum_{\alpha} t_{\alpha}
\sum_{\beta + \gamma=\alpha} N_{\beta \gamma} \omega^{\beta} \wedge \omega^{\gamma}
\pmod{\text{semibasic terms}}.
\end{gather*}
In particular, for any noncompact negative root $\alpha$,
\[
\LieDer_{X_{-\alpha}} t_{\alpha} = 2 \frac{\left<\delta,\alpha\right>}{\left<\alpha,\alpha\right>}.
\]
Since $-\alpha$ is a positive root, the corresponding root vector
lies in $\mathfrak{p}$. Moreover, this root vector lies in the
nilpotent radical of $\mathfrak{p}$, since it is a positive root.
The nilpotent radical acts trivially on~$\Lambda^{\text{top},0}\left(\mathfrak{g}/\mathfrak{p}\right)$,
as nilpotent groups have no nontrivial 1-dimensional representations.
Every vector in $\mathfrak{p}$ gives rise to a vector f\/ield
giving the associated inf\/initesimal action, and for the root
vector of $-\alpha$, this vector f\/ield is $X_{-\alpha}$,
by def\/inition of a Cartan geometry. Since the root
vector lies in the Lie algebra of $P_0$, the vector f\/ield
$X_{-\alpha}$ generates a 1-parameter subgroup of $P_0$.
The vector f\/ield $X_{-\alpha}$ is therefore
tangent to the f\/ibers of $E' \to M$.

On the f\/ibers the vector f\/ield $X_{-\alpha}$ is a left invariant vector f\/ield.
Therefore $X_{-\alpha}$ is complete. Starting at any point of $E'$,
we can move in the direction $X_{-\alpha}$ of the nilpotent part
of the structure group, altering the value of $t_{\alpha}$ at a constant rate
until it reaches $0$. Indeed $t_{\alpha}$ is acted on by the nilradical
of the structure group as translations in the left action on $E'$. The set of points
$E'' \subset E'$ on which all $t_{\alpha}$ vanish is a smooth embedded
submanifold, because its tangent space is cut out by equations
\[
\frac{\left<\delta,\alpha\right>}{\left<\alpha,\alpha\right>} \omega^{-\alpha} = \text{semibasic},
\]
for all noncompact negative roots $\alpha.$ The structure group is
reduced to
a reductive algebraic group, since we have eliminated the nilradical
of the original structure group, leaving only the root spaces
$\alpha$ for which neither $\alpha$ nor $-\alpha$ is
noncompact negative, i.e.
the root spaces of the maximal reductive subgroup $MA$
of the structure group $P$. We have also eliminated the
part of~$A$ which acts nontrivially on the holomorphic
volume forms, so our structure group is now~$MA^0$,
with~$A^0$ the subgroup of~$A$ f\/ixing a volume form on
$\mathfrak{g}/\mathfrak{p}$.
\end{proof}

\begin{remark}
On a complex manifold with a meromorphic section of the canonical
bundle, it would be interesting to consider what happens to this
argument as we approach the zeroes or poles of the meromorphic
section.
\end{remark}

\begin{corollary}\label{corollary:Weyl}
If a complex manifold bears a holomorphic parabolic geometry
and a holomorphic volume form, then it admits a canonical
holomorphic Weyl structure.
\end{corollary}
\begin{proof}
Suppose that $E \to M$ is a holomorphic parabolic
geometry modelled on $G/P$. Write the Langlands
decomposition of $P$ as $P=MAN$.
Let $G^0=P/N$.
Recall that a holomorphic \emph{Weyl structure} is a
$G^0$-equivariant holomorphic section of the bundle
$E \to E/N$ \cite{Cap/Slovak:2009}. By Theorem~\ref{theorem:Reduction},
we have a principal right $MA^0$-subbundle $E'' \subset E$,
which induces a principal right $MA$-subbundle
$E'' \times_{MA^0} MA \subset E$.
Clearly $P/N=MA$ and the map $E_0 \to E/N$ is an
isomorphism of $MA$-bundles.
\end{proof}

\section{Parabolic geometries on tori}

\begin{example}\label{example:TranslationInvariant}
Suppose that $L$ is a Lie group with Lie algebra $\LieL$,
and write the left invariant Maurer--Cartan 1-form
on $L$ as $\ell^{-1}   d\ell$.
Suppose that $G$ is a Lie group and $H \subset G$
is a closed subgroup. Take any linear injection $t : \LieL \to \LieG$
so that the image is complementary to $\LieH$.
Then let $M=L$ and $E=M \times H$. Let $\omega \in \nForms{1}{E} \otimes  \LieG$
be the 1-form
\[
\omega = h^{-1}   dh + \Ad(h)^{-1}\big(t   \ell^{-1}   d\ell\big).
\]
It is easy to check that $\omega$ is the Cartan
connection of Cartan geometry on $L$ modelled on $G/H$.
The group $L$ acts as Cartan geometry automorphisms
by the obvious action. The curvature is
\[
d \omega + \frac{1}{2}\left[\omega,\omega\right]
=
\frac{1}{2}
\Ad(h)^{-1}\left(
\left[
t \ell^{-1}   d\ell,
t \ell^{-1}   d\ell
\right]
-
t
\left[
\ell^{-1}   d\ell,
\ell^{-1}   d\ell
\right]
\right).
\]
In particular, the Cartan geometry is f\/lat if and only if
$t$ is a Lie algebra homomorphism. See~\cite{Hammerl:2007}
for details. Clearly the group $\Aut(L) \times \Aut(G,H)$
acts as isomorphisms of these geometries.
\end{example}

\begin{lemma} Every left invariant Cartan geometry on a
Lie group is equivariantly isomorphic to one
constructed as in Example~{\rm \ref{example:TranslationInvariant}}.
\end{lemma}
\begin{proof}
Take any Cartan geometry $E \to L$ on a Lie group $L$.
Suppose that $L$ acts on $E$ lifting its left action
on itself, commuting with the $H$-action,
and preserving the Cartan connection.
Pick any point $e_0 \in E$. Map
\[
\left(\ell,h\right) \in L \times H \mapsto \ell e_0 h.
\]
Clearly this is an isomorphism of principal bundles,
so from now on we take $E=L \times H$.
Unwinding the def\/inition of a Cartan connection immediately yields
that every Cartan connection on the bundle $L \times H \to L$ has the form
\[
\omega = h^{-1}   dh + \Ad(h)^{-1} \left(\gamma\right),
\]
where $\gamma$ is a 1-form on $L$ valued in $\LieG$,
so that $\gamma+\LieH$ is a linear isomorphism,
i.e.\ $\gamma_{\ell} : T_{\ell} L \to \gtot$ is a
linear injection, for each $\ell \in L$, and
$T_{\ell} L \to \gtot \to \gtot/\gstruc$ is a linear isomorphism.
By translation invariance
\[
\gamma = t \ell^{-1}   d\ell,
\]
for a unique linear map $t$ and
$t \to \gtot \to \gtot/\gstruc$ is a linear isomorphism.
\end{proof}

\begin{theorem}\label{theorem:translation}
Every holomorphic parabolic geometry on any complex torus is translation
inva\-riant, and obtained by the construction of
Example~{\rm \ref{example:TranslationInvariant}}
$($where the Lie group $L$ is the complex torus itself$)$.
More generally, if $M$ is any compact complex
manifold with holomorphically trivial tangent bundle,
then $M=L/\Gamma$ for some discrete
subgroup $\Gamma \subset L$ of a complex
Lie group~$L$. Every holomorphic parabolic geometry
on $M$ is obtained by the construction of
Example~{\rm \ref{example:TranslationInvariant}}
applied to~$L$, and then quotiented by~$\Gamma$.
\end{theorem}

\begin{proof}
If $M$ is a compact complex manifold with holomorphically
trivial tangent bundle, then any basis of the
holomorphic tangent space $T_m M$ at any point $m \in M$
extends to a framing by holomorphic vector f\/ields, say $X_1, X_2, \dots, X_n$.
These then generate a transitive complex Lie group action,
say of a complex Lie group $L$. Brackets of these
holomorphic vector f\/ields can be rewritten in terms
of the vector f\/ields themselves
\[
\left[X_i,X_j\right]=\sum c^k_{ij} X_k,
\]
since the $X_k$ form a basis. The holomorphic functions
$c^k_{ij}$ are constant because $M$ is compact. Therefore $L$ has the same
dimension as $M$, and so $M = L/\Gamma$ \cite{Wang:1954}.

Following the proof of Theorem~\ref{theorem:Reduction}, the structure group of any
holomorphic parabolic geometry $E \to M$ reduces to a reductive group, $G''$
on some subbundle $E'' \subset E$. The Cartan connection~$\omega$ splits into a sum corresponding to the splitting
of $\mathfrak{g}$ into $G''$-invariant subspaces,
and~$\omega''$ (the part valued in $\mathfrak{g}''$)
is a connection form for $E'' \to M$. Take a global
holomorphic coframing on~$M$, i.e.\ a set of linearly independent 1-forms
$\xi^{\alpha}$ forming a basis of each cotangent space of~$M$, for
$\alpha$ varying over noncompact negative roots.
Def\/ine a map $e \in E'' \to h
\in \GL{\gtot/\gstruc}$, by $\omega^{\alpha}=h^{\alpha}_{\beta}
\xi^{\beta}$ (for $\alpha$ and $\beta$ varying over
noncompact negative roots),
and $h(e) =  (h^{\alpha}_{\beta} )$ in the basis~$X_{\alpha}$ for the sum of noncompact negative
root spaces. Under right
$G''$-action,
\[
h\left(r_{g}e\right)=g^{-1}h(e)
\]
for $g \in G''.$ Therefore the quotient map $E \!\to\! \GL{\gtot/\gstruc}/G''$
descends to a map $M \!\to\! \GL{\gtot/\gstruc}/G''$.
The quotient $\GL{\gtot/\gstruc}/G''$ is an af\/f\/ine variety:
see Mumford et al.~\cite[Theorem~1.1, p.~27]{Mumford/Fogarty/Kirwan:1994}
and Procesi~\cite[Theorem~2, p.~556]{Procesi:2007}.
Af\/f\/ine coordinate functions will pull back to functions on $M$,
and therefore must be constant.
Therefore the map $M \to \GL{\gtot/\gstruc}/G''$ is constant.
We have an isomorphism
\[
e \in E'' \to \left(\pi(e),h(e)\right) \in M \times G'',
\]
trivializing the bundle $E''$. We can therefore assume
that $E'' = M \times G''$ and $E = M \times H$.
Again unwinding the def\/inition of a Cartan
geometry,
\[
\omega = h^{-1}   dh + \Ad(h)^{-1}\gamma,
\]
and $\gamma \in \nForms{1}{M} \otimes \LieG$.
The functions $X_i \hook \gamma$
are holomorphic functions on $M$,
so constant, so $\gamma$ pulls back to $L$ to be
$\gamma=t \ell^{-1}   d\ell$ for some constant
linear map $t : \LieL \to \LieG$.
\end{proof}

\begin{remark}
This theorem is more remarkable if one remembers that there
are holomorphic Cartan geometries on complex tori which are
not translation invariant.
\end{remark}

\begin{remark}
Suppose that $G$ is a complex semisimple Lie group and $P \subset G$
is a parabolic subgroup. If $\LieG$ has no Abelian subalgebra complementary
to $\LieP$, then there is no f\/lat holomorphic $G/P$-geometry
on any complex torus. Indeed this occurs for all $G/P$
which are not compact Hermitian symmetric spaces.
More generally, for any $G$ and $P$, there is some complex linear subspace
in $\LieG$ complementary to $\LieP$ which is not an Abelian
subalgebra. Therefore there is a~holomorphic $G/P$-geometry
on any complex torus of the appropriate dimension which
is not f\/lat.
\end{remark}

\begin{corollary}
If a compact K\"ahler manifold with $c_1=0$ bears a parabolic
geometry, then it is covered by a torus, and the parabolic geometry
pulls back to a translation invariant parabolic geometry on the
torus.
\end{corollary}

\begin{definition}\label{definition:lift}
Suppose that $\Gstruc_- \subset \Gstruc_+$ are two closed subgroups
of a Lie group $\Gtot$, so that we have a f\/iber bundle map
$\Gtot/\Gstruc_- \to \Gtot/\Gstruc_+$. Let $E \to M$ be a Cartan geometry modelled on~$\Gtot/\Gstruc_+$. Then $E \to E/\Gstruc_-$ is a Cartan geometry
modelled on $\Gtot/\Gstruc_-,$ called the \emph{lift} of the Cartan
geometry on $M$.
\end{definition}

\begin{corollary}
On any compact homogeneous K\"ahler manifold, all parabolic
geometries are lifted $($as in Definition~{\rm \ref{definition:lift})} from a
translation invariant geometry on a torus $($constructed as in
Example~{\rm \ref{example:TranslationInvariant})}. In particular, all
such parabolic geometries are homogeneous.
\end{corollary}

\begin{proof}
Borel and Remmert \cite{Borel/Remmert:1961} proved that every compact
homogeneous K\"ahler manifold is a product of a torus and a rational
homogeneous variety. The rational homogeneous variety bears rational
curves just when it has positive dimension. These rational curves
ensure that the parabolic geometry is lifted from lower dimension
(see Biswas and McKay \cite{Biswas/McKay:2010c}), quotienting out the rational
homogeneous variety entirely.
\end{proof}

\begin{remark}
Any parabolic geometry on any rational homogeneous variety
(or more generally, on any compact rationally connected complex
manifold) is f\/lat and isomorphic to its model
(see Biswas and McKay~\cite{Biswas/McKay:2010c}).
\end{remark}

\section{Conclusion}

The classif\/ication of holomorphic parabolic geometries
on compact complex manifolds with \mbox{$c_1>0$} is complete \cite{Biswas/McKay:2010c}
(and more generally the classif\/ication on rationally connected
compact complex manifolds).
Above we give the classif\/ication of holomorphic parabolic geometries
on compact K{\"a}hler manifolds with $c_1=0$. The classif\/ication for
$c_1<0$ must be more dif\/f\/icult, as there are many locally
symmetric varieties and many holomorphic projective connections
on compact Riemann surfaces of genus $g \ge 2$.

We propose a conjecture:
if $M$ is a compact complex manifold with $c_1<0,$ then either (1)~$M$~admits no parabolic geometry, or (2)~$M$~admits a parabolic
geometry modelled on a~compact Hermitian symmetric space
$\Gtot/\Gstruc$ and $M$ is covered by the
noncompact dual of that symmetric space. In case~(2), every
parabolic geometry on $M$ modelled on $\Gtot/\Gstruc$ is f\/lat,
and if the factorization of $G/P$ into a product of irreducible
compact Hermitian symmetric spaces has no factor of dimension~1
(i.e.\ isomorphic to $\Proj{1}$) then the parabolic geometry
on $M$ pulls back to the standard f\/lat parabolic geometry on
the noncompact dual.

The above conjecture was proven by Kobayashi \cite{Kobayashi/Horst:1983}
with the additional hypothesis
that the parabolic geometry is a normal conformal connection or a
normal projective connection.

There are examples of smooth complex projective varieties
which are not locally symmetric and which have holomorphic
projective connections~\cite{Jahnke/Radloff:2002}.
These examples stand in the way of any obvious conjecture
as to which smooth complex projective varieties have
holomorphic parabolic geometries. They have $c_1 \le 0$
but not $c_1 < 0$ or $c_1 = 0$. It might help
to be able to limit the possible models as follows.

We propose another conjecture:
suppose that $M$ is a compact connected K{\"a}hler manifold,
bearing a holomorphic parabolic geometry. Then
the canonical bundle of $M$ is not pseudoef\/fective
if and only if the parabolic geometry drops
(in the sense of~\cite{Biswas/McKay:2010c}) to a
holomorphic parabolic geometry on a lower
dimensional compact K{\"a}hler manifold.
(We can essentially ignore such parabolic
geometries in any classif\/ication.)
In the other direction, if the canonical bundle
of~$M$ is pseudoef\/fective, and the parabolic
geometry is regular at at least one point,
then the model is a compact Hermitian symmetric
space. (Consequently we hope to reduce the
classif\/ication to the classif\/ication of holomorphic
parabolic geometries modelled on compact Hermitian
symmetric spaces, a topic which presumably lies close
to the study of locally Hermitian symmetric varieties.)

We propose another conjecture, concerning the complex torus:
on any complex torus, a~holomorphic Cartan geometry
is translation invariant if and only if it is locally homogeneous.

We expect that foliations play a fundamental role in
the phenomenon of translation inva\-riance.
We conjecture:
suppose that $G$ is a complex Lie group and $H \subset G$
is a closed complex subgroup.
Then either (1) every $G$-invariant holomorphic foliation
of $G/H$ has a $G$-invariant holomorphic complementary subbundle
of the tangent bundle
and every holomorphic Cartan geometry modelled on $G/H$
on any complex torus is translation invariant
or (2) there is a~$G$-invariant holomorphic foliation of $G/H$,
with no invariant holomorphic complementary subbundle of the tangent
bundle, and there is a complex Abelian variety $A$ and
a holomorphic Cartan geometry on $A$ modelled on $G/H$
which is not translation invariant.

\subsection*{Acknowledgements}
This material is based upon works supported by the Science Foundation
Ireland under Grant No.~MATF634.

\pdfbookmark[1]{References}{ref}
 \LastPageEnding


\begin{thebibliography}{99}

\footnotesize\itemsep=0pt

\bibitem{Atiyah:1957}
 Atiyah M.F.,
Complex analytic connections in f\/ibre bundles,
\href{http://dx.doi.org/10.2307/1992969}{{\em Trans. Amer. Math. Soc.}} {\bf 85} (1957), 181--207.

\bibitem{Lascoux/Berger:1970}
  Berger M., Lascoux A.,
Vari\'et\'es K\"ahleriennes compactes,
{\it Lecture Notes in Mathematics}, Vol.~154, Springer-Verlag, Berlin, 1970.

\bibitem{Biswas/McKay:2010a}
 Biswas I.,  McKay B.,
 Holomorphic Cartan geometries and Calabi--Yau manifolds,
\href{http://dx.doi.org/10.1016/j.geomphys.2009.12.011}{{\em J.~Geom. Phys.}} {\bf 60} (2010), 661--663,
\href{http://arxiv.org/abs/0812.3978}{arXiv:0812.3978}.

\bibitem{Biswas/McKay:2010c}
Biswas I.,  McKay B.,
Holomorphic Cartan geometries and rational curves,
\href{http://arxiv.org/abs/1005.1472}{arXiv:1005.1472}.

\bibitem{Biswas/McKay:2010b}
Biswas I.,  McKay B.,
Holomorphic Cartan geometries, Calabi--Yau manifolds and rational curves,
\href{http://dx.doi.org/10.1016/j.difgeo.2009.09.003}{{\em Differential Geom. Appl.}} {\bf 28} (2010), 102--106,
\href{http://arxiv.org/abs/1009.5801}{arXiv:1009.5801}.

\bibitem{Borel/Remmert:1961}
Borel A., Remmert R.,
\"{U}ber kompakte homogene K\"ahlersche Mannigfaltigkeiten,
\href{http://dx.doi.org/10.1007/BF01471087}{{\em Math. Ann.}} {\bf 145}  (1961/1962), 429--439.

\bibitem{Cap/Slovak:2009}
\v{C}ap A., Slov{\'a}k J.,
Parabolic geometries.~I.~Background and general theory,
{\em Mathematical Surveys and Monographs}, Vol.~154,
American Mathematical Society, Providence, RI, 2009.


\bibitem{Dumitrescu:2001}
Dumitrescu S.,
Structures g\'eom\'etriques holomorphes sur les vari\'et\'es   complexes compactes,
\href{http://dx.doi.org/10.1016/S0012-9593(01)01070-9}{{\em Ann. Sci. \'Ecole Norm. Sup.~(4)}} {\bf 34} (2001), 557--571.

\bibitem{Dumitrescu:2010b}
Dumitrescu S.,
Connexions af\/f\/ines et projectives sur les surfaces complexes   compactes,
\href{http://dx.doi.org/10.1007/s00209-008-0465-8}{{\em Math.~Z.}} {\bf 264} (2010), 301--316,
\href{http://arxiv.org/abs/0805.2816}{arXiv:0805.2816}.

\bibitem{Dumitrescu:2010}
Dumitrescu S.,
Killing f\/ields of holomorphic Cartan geometries,
\href{http://dx.doi.org/10.1007/s00605-009-0135-x}{{\em Monatsh. Math.}} {\bf 161} (2010), 145--154,
\href{http://arxiv.org/abs/0902.2193}{arXiv:0902.2193}.

\bibitem{Fulton/Harris:1991}
 Fulton W.,  Harris J.,
  Representation theory. A~f\/irst course,
  {\em Graduate Texts in   Mathematics}, Vol.~129,
 Springer-Verlag, New York, 1991.

\bibitem{Gunning:1978}
Gunning R.C.,
On uniformization of complex manifolds: the role of connections,
{\em Mathematical Notes}, Vol.~22,
Princeton University Press, Princeton, N.J., 1978.

\bibitem{Hammerl:2007}
Hammerl M.,
Homogeneous Cartan geometries,
{\em Arch. Math. (Brno)} {\bf 43} (2007), 431--442,
\href{http://arxiv.org/abs/math.DG/0703627}{math.DG/0703627}.

\bibitem{Igusa:1954}
 Igusa J.-I.,
 On the structure of a certain class of Kaehler varieties,
\href{http://dx.doi.org/10.2307/2372709}{{\em Amer.~J. Math.}}  {\bf 76} (1954), 669--678.

\bibitem{Inoue/Kobayashi/Ochiai:1980}
  Inoue M.,  Kobayashi S.,  Ochiai T.,
 Holomorphic af\/f\/ine connections on compact complex surfaces,
 {\em J.~Fac. Sci. Univ. Tokyo Sect. IA Math.} {\bf 27} (1980), 247--264.

\bibitem{Jahnke/Radloff:2002}
 Jahnke P., Radlof\/f I.,
 Threefolds with holomorphic normal projective connections,
\href{http://dx.doi.org/10.1007/s00208-004-0406-8}{{\em Math. Ann.}} {\bf 329} (2004), 379--400,
\href{http://arxiv.org/abs/math.AG/0210117}{math.AG/0210117}.

\bibitem{Jahnke/Radloff:2004}
 Jahnke P., Radlof\/f I.,
Projective threefolds with holomorphic conformal structure,
\href{http://dx.doi.org/10.1142/S0129167X05003028}{{\it Internat.~J. Math.}} {\bf 16} (2005), 595--607,
\href{http://arxiv.org/abs/math.AG/0406113}{math.AG/0406113}.

\bibitem{Klingler:1998}
Klingler B.,
Structures af\/f\/ines et projectives sur les surfaces complexes,
{\em Ann. Inst. Fourier (Grenoble)} {\bf 48} (1998), 441--477.

\bibitem{Klingler:2001}
Klingler B.,
Un th\'eor\`eme de rigidit\'e non-m\'etrique pour les vari\'et\'es localement sym\'etriques hermitiennes,
\href{http://dx.doi.org/10.1007/s00014-001-8320-0}{{\em Comment. Math. Helv.}} {\bf 76} (2001), 200--217.

\bibitem{Knapp:2002}
 Knapp A.W.,
Lie groups beyond an introduction, 2nd ed.,
 {\em Progress   in Mathematics}, Vol.~140,
  Birkh\"auser Boston Inc., Boston, MA, 2002.

\bibitem{Kobayashi/Horst:1983}
  Kobayashi S.,  Horst C.,
 Topics in complex dif\/ferential geometry,
in  Complex Dif\/ferential Geometry,   {\em DMV Sem.}, Vol.~3,
 Birkh\"auser, Basel, 1983, 4--66.

\bibitem{KobayashiOchiai:1980}
 Kobayashi S., Ochiai T.,
Holomorphic projective structures on compact complex surfaces,
\href{http://dx.doi.org/10.1007/BF01387081}{{\em Math. Ann.}} {\bf 249} (1980), 75--94.

\bibitem{Kobayashi/Ochiai:1981ii}
Kobayashi S., Ochiai T.,
 Holomorphic projective structures on compact complex surfaces.~II,
\href{http://dx.doi.org/10.1007/BF01451931}{{\em Math. Ann.}} {\bf 255} (1981), 519--521.

\bibitem{Kobayashi/Ochiai:1981}
Kobayashi S., Ochiai T.,
 Holomorphic structures modeled after compact Hermitian symmetric spaces,
in Manifolds and Lie Groups (Notre Dame, Ind., 1980),
 {\em Progr. Math.}, Vol.~14,  Birkh\"auser, Boston, Mass., 1981, 207--222.

\bibitem{Kobayashi/Ochiai:1982}
Kobayashi S., Ochiai T.,
 Holomorphic structures modeled after hyperquadrics,
\href{http://dx.doi.org/10.2748/tmj/1178229159}{{\em T\^ohoku Math.~J.~(2)}} {\bf 34} (1982), 587--629.

\bibitem{McKay:2011c}
 McKay B.,
 Characteristic forms of complex Cartan geometries,
\href{http://dx.doi.org/10.1515/ADVGEOM.2010.044}{{\em Adv. Geom.}} {\bf 11} (2011), 139--168,
\href{http://arxiv.org/abs/0704.2555}{arXiv:0704.2555}.

\bibitem{McKay:2011}
McKay B.,
Holomorphic Cartan geometries on uniruled surfaces,
\href{http://dx.doi.org/10.1016/j.crma.2011.07.021}{{\em C. R. Acad. Sci. Paris}} {\bf 349} (2011), 893--896,
\href{http://arxiv.org/abs/1105.4732}{arXiv:1105.4732}.

\bibitem{McKay/Pokrovskiy:2010}
McKay B.,   Pokrovskiy A.,
 Locally homogeneous geometric structures on Hopf surfaces,
\href{http://dx.doi.org/10.1512/iumj.2010.59.4178}{{\em Indiana Univ. Math.~J.}} {\bf 59} (2010), 1491--1540,
\href{http://arxiv.org/abs/0910.0369}{arXiv:0910.0369}.

\bibitem{Mumford/Fogarty/Kirwan:1994}
Mumford D., Fogarty J.,  Kirwan F.,
 Geometric invariant theory,  3rd ed., {\em Ergebnisse der
  Mathematik und ihrer Grenzgebiete (2)}, Vol.~34,
 Springer-Verlag, Berlin,   1994.

\bibitem{Procesi:2007}
 Procesi C.,
Lie groups. An approach through invariants and representations, {\it Universitext}, Springer, New York, 2007.

\bibitem{Serre:2001}
 Serre J.-P.,
 Complex semisimple Lie algebras,
{\it Springer Monographs in Mathematics}, Springer-Verlag, Berlin, 2001.


\bibitem{Sharpe:1997}
 Sharpe R.W.,
Dif\/ferential geometry. Cartan's generalization of Klein's Erlangen program,
{\em Graduate Texts in Mathematics}, Vol.~166,
 Springer-Verlag, New York, 1997.

\bibitem{Wang:1954}
 Wang H.-C.,
 Closed manifolds with homogeneous complex structure,
\href{http://dx.doi.org/10.2307/2372397}{{\em Amer.~J. Math.}} {\bf 76} (1954), 1--32.

\bibitem{Yau:1977}
 Yau S.T.,
 Calabi's conjecture and some new results in algebraic geometry,
 {\em Proc. Nat. Acad. Sci. U.S.A.} {\bf 74} (1977), 1798--1799.

\end{thebibliography}
\end{document}